\documentclass[12pt]{amsart}
\usepackage{amsfonts, amssymb, amsmath, amscd, bbold, graphics}

\newcommand{\C}{\mathbb{C}}
\newcommand{\tr}[1]{\mathrm{tr}(#1)}
\newcommand{\xb}{\mathbf{x}}
\newcommand{\yb}{\mathbf{y}}
\newcommand{\zb}{\mathbf{z}}
\newcommand{\pol}{\mathrm{pol}}
\newcommand{\ub}{\mathbf{u}}
\newcommand{\vb}{\mathbf{v}}
\newcommand{\wb}{\mathbf{w}}
\newcommand{\dg}[1]{|\!|#1|\!|}
\newcommand{\tdg}[1]{|\!|#1|\!|_{\mathrm{tr}}}
\newcommand{\aq}{/\!\!/}
\newcommand{\X}{\mathfrak{X}}
\newcommand{\R}{\mathfrak{R}}
\newcommand{\hm}{\mathrm{Hom}}
\newcommand{\SL}{\mathrm{SL}(3,\C)}
\newcommand{\ot}{\mathrm{Out}}
\newcommand{\F}{\mathtt{F}}
\newcommand{\G}{\mathfrak{G}}
\newcommand{\xt}{\mathtt{x}}
\newcommand{\M}{\mathbf{M}}
\newcommand{\id}{\mathbb{I}}
\newcommand{\ti}[1]{t_{(#1)}}
\newcommand{\dP}[1]{\frac{\partial P}{\partial \ti{#1}}}
\newcommand{\dQ}[1]{\frac{\partial Q}{\partial \ti{#1}}}

\newtheorem{thm}{Theorem}
\newtheorem{lem}[thm]{Lemma}
\newtheorem{cor}[thm]{Corollary}
\newtheorem{rem}[thm]{Remark}
\newtheorem{prop}[thm]{Proposition}
\newtheorem{fact}[thm]{Fact}

\title{Generators, Relations and Symmetries in Pairs of $3\times 3$ Unimodular Matrices}
\author{Sean Lawton}
\address{Mathematics Department, Kansas State University, Manhattan, KS 66506}
\address{E-mail address: $\mathtt{ slawton@math.ksu.edu}$}
\address{URL: $\mathtt{http://www.math.ksu.edu/\sim slawton}$}
\date{January 22, 2007}

\begin{document}

\maketitle

\begin{abstract}
Denote the free group on two letters by $\F_2$ and the $\SL$-representation variety of $\F_2$ by $\R=\hm(\F_2,\SL)$.  
There is a $\SL$-action on the coordinate ring of $\R$, and the geometric points of the subring of invariants is an affine 
variety $\X$.  We determine explicit minimal generators and defining relations for the subring of invariants and show $\X$ is a degree $6$ hyper-surface 
in $\C^9$ mapping onto $\C^8$.  Our choice of generators exhibit $\ot(\F_2)$ symmetries which allow for a succinct expression of the defining relations.
\end{abstract}

\section{Introduction}
The purpose of this paper is to describe a minimal generating set and defining relations for the 
ring of invariants $$\C[\SL\times\SL]^{\SL}.$$  This generating set exhibits symmetries which allow for an explicit and succinct expression of the invariant ring 
as a quotient.  

Explicit minimal generators have been found by \cite{We} and graphically by \cite{Si}; in an unpublished calculation \cite{PX} independently describe the defining relations.   Our 
treatment provides the most succinct and transparent description by uncovering symmetries which provide a framework for generalization.

A related algebra, however different, is the ring of invariants of pairs of $3\times 3$ matrices $M_3(\C)\times M_3(\C)$ under simultaneous conjugation.  The algebra of invariants 
$\C[M_3(\C)\times M_3(\C)]^{\mathrm{GL}(3,\C)}$ comes to bear on the algebra of invariants $\C[\SL\times\SL]^{\SL}$ by restriction.  On the other hand, the former ring of invariants may 
be described, in part, by $\C[\frak{sl}(3)\times \frak{sl}(3)]^{\SL}$; the infinitesimal invariants in the latter ring.  In this more general context, similar 
questions about generators and relations have been addressed.  In particular, explicit minimal generators were first found by \cite{Du} in $1935$, and 
later by \cite{SR, MS, T}.  The much more general results of \cite{AP} additionally provide minimal generators.  However, \cite{N}, and later \cite{ADS} were the first to explicitly describe 
the defining relations.  For the state-of-the-art, see \cite{DF}.

We now describe the main results of this paper.  Let $\X$ be the variety whose coordinate ring is $\C[\X]=\C[\SL\times\SL]^{\SL}$.  Theorem $\ref{ranktwo}$ 
asserts that $\X$ is isomorphic to a degree $6$ affine hyper-surface in $\C^9$ which generically maps $2$-to-$1$ onto $\C^8$.  Next, Theorem 
$\ref{singulartheorem}$ explicitly describes the singular locus of $\X$, and examples of non-singular representations in the branching 
locus are constructed.  Lastly, Theorem $\ref{symmetrytheorem}$ describes an $8$-fold symmetry on $\C[\X]$ which at once characterizes the algebraically 
independent generators and allows for a surprisingly simple description of the defining relations.

We hope that this paper will be of interest to algebraic-geometers, ring theorists, and geometers alike.  In particular, results in this paper have recently been used in work concerning 
the hyperbolic geometry of sphereical CR manifolds (see \cite{Sc}).  With this in mind, some of the exposition, for instance, may be ``well-known'' to a ring theorist but perhaps not 
to an algebraic-geometer or a geometer.  The reader is encouraged to skip such exposition, as appropriate.

\section*{Acknowledgements}
The author thanks William Goldman for overall guidance, substantial suggestions, and edits that always improved the quality of earlier drafts of this paper.
He thanks Richard Schwartz, John Millson, and the University of Maryland's VIGRE-NSF program for generously supporting this work.
In particular, he thanks John Millson for suggesting an analysis of the singular locus, which has greatly improved the quality of this paper.  The author has benefited from fruitful
conversations with Ben Howard, Elisha Peterson, Adam Sikora, and Joseph Previte, and thanks them for their time and insight.  He especially thanks Joseph Previte and Eugene Xia for
generously sharing their calculations with him.  Lastly, the author thanks the reviewer, whose suggestions and corrections have helped make this paper more readable.

\section{$\SL$ Invariants}
\subsection{Algebraic Structure of $\SL$}

The group $\SL$ has the structure of an algebraic set since it is the zero set of the polynomial $$D=\det\left(
\begin{array}{ccc}
x_{11} & x_{12} & x_{13}\\
x_{21} & x_{22} & x_{23}\\
x_{31} & x_{32} & x_{33}\\
\end{array}\right)-1$$ on $\mathbb{C}^9$.  Here $x_{ij}\in
\C[x_{11},x_{12},x_{13},x_{21},x_{22},x_{23},x_{31},x_{32},x_{33}]$, the polynomial ring over $\C$ in 
$9$ indeterminates.  As such denote $\SL$ by $\G$.  The coordinate ring of $\G$ is given by $$\mathbb{C}[\G]=\mathbb{C}[x_{ij}\ | \ 
1\leq i,j \leq 3]/(D).$$  Since $D$ is irreducible, $(D)$ is a prime ideal.  So the algebraic set $\G$ is in fact an affine variety.  

\subsection{Representation and Character Varieties of a Free Group}\label{repchar}

Let $\F_r$ be the free group of rank $r$ generated by $\{\xt_1,...,\xt_r\}.$  The map $$\hm(\F_r,\G)\longrightarrow
\G^{\times r}$$ defined by sending $$\rho \mapsto (\rho(\xt_1),\rho(\xt_2),...,\rho(\xt_r))$$ is a
bijection.  Since $\G^{\times r}$ is the $r$-fold product of irreducible algebraic sets, $\G^{\times r} \cong \hm(\F_r,\G)$ is an affine 
variety.
 
As such $\hm(\F_r, \G)$ is denoted by $\R$ and referred to as the $\SL$-{\it representation variety of $\F_r$}.  

Let $\C[\R]$ be the coordinate ring of $\R$.  Our preceding remarks imply 
$\C[\R]\cong \C[\G]^{\otimes r}$.  
For $1\leq k\leq r,$ define  a {\it generic matrix} of the complex polynomial ring in $9r$ indeterminates by
$$\xb_k=\left(
\begin{array}{ccc}
x^k_{11} & x^k_{12} & x^k_{13}\\
x^k_{21} & x^k_{22} & x^k_{23}\\
x^k_{31} & x^k_{32} & x^k_{33}\\
\end{array}\right).$$   
Let $\Delta$ be the ideal $(\det(\xb_k)-1\ |\ 1\leq k\leq r)$ in $\C[\R]$.  Then 
$$\C[\R]=\mathbb{C}[x^k_{ij}\ | \ 1\leq i,j \leq 3,\ 1\leq k\leq r]/\Delta.$$

Let $(\xb_1,\xb_2,...,\xb_r)$ be an $r$-tuple of generic matrices.  An element $f\in \C[\R]$ is a function defined in 
terms of such $r$-tuples.  There is a polynomial $\G$-action on $\C[\R]$ given by diagonal conjugation; that is, for $g\in \G$ $$g\cdot 
f(\xb_1,\xb_2,...,\xb_r)= f(g^{-1}\xb_1 g,...,g^{-1}\xb_r g).$$  
The subring of invariants of this action $\C[\R]^{\G}$ is a finitely 
generated $\C$-algebra (see \cite{D,P,R}).  Consequently, the {\it character variety} 
$$\X=\mathrm{Spec}_{max}(\C[\R]^{\G})$$ is the irreducible algebraic set whose coordinate ring is the 
ring of invariants.   For $r>1$, the Krull dimension of $\X$ is $8r-8$ since generic elements have zero
dimensional isotropy (see \cite{D}, page $98$).  More generally, the dimension of the ring of invariants $\C[M_n(\C)^{\times r}]^{\mathrm{GL}(n,\C)}$ is $n^2(r-1)+1$ (see \cite{DF}).  
Consequently, the dimension of $\C[\frak{sl}(n)^{\times r}]^{\mathrm{SL}(n,\C)}$, which equals that of $\C[\X]$, is $(n^2-1)(r-1)$. 
  
There is a regular map $\R\stackrel{\pi}{\to}\X$ which factors through $\R/\G$:  let $\frak{m}$ be a maximal ideal
corresponding to a point in $\R$, then the composite isomorphism $\C\to \C[\R]\to \C[\R]/\frak{m}$ implies that the composite map $\C\to
\C[\R]^\G\to \C[\R]^\G/(\frak{m}\cap \C[\R]^\G)$ is an isomorphism as well.  Hence the contraction $\frak{m}\cap \C[\R]^\G$ is maximal, and since for any 
$g\in \G$, $\left(g\frak{m}g^{-1}\right)\cap \C[\R]^\G=\frak{m}\cap \C[\R]^\G$, $\pi$ factors through $\R/\G$ (see \cite{E}, page $38$).
Although $\R/\G$ is not generally an algebraic set, $\X$ is the categorical quotient $\R\aq\G$, and since $\G$ is a (geometrically) reductive algebraic 
group $\pi$ is surjective, maps closed $\G$-invariant sets to closed sets, and separates distinct closed orbits (see \cite{D}).

\section{Trace Identities for Matrices}

Let $\F^{+}_r$ be the free monoid generated by $\{\xt_1,...,\xt_r\}$, and let $\M^+_r$ be the monoid  
generated by $\{\xb_1,\xb_2,...,\xb_r\},$ as defined in section $\ref{repchar}$, under matrix multiplication and with identity $\id$ the $3\times 3$ 
identity matrix.  There is a surjection $\F^{+}_r\to \M^+_r$, defined by mapping $\xt_i\mapsto \xb_i$.   
Let $\mathbf{w} \in \M^+_r$ be the image of $\mathtt{w}\in \F^{+}_r$ under this map.  Further, let $|\cdot|$ be the function 
that takes a reduced word in $\F_r$ to its word length.  Then by \cite{P,R}, we know $\C[\X]$ is not only finitely generated, but in fact 
generated by 
\begin{equation}\label{generators}
\{\tr{\mathbf{w}}\ |\ \mathtt{w}\in \F^{+}_r, \ |\mathtt{w}|\leq 7\}.  
\end{equation}

More generally, the length of the generators is bounded by the class of nilpotency of nil algebras of class $n$.  With respect to matrix algebras, $n$ is the size of the matrices 
under consideration.  The best known upper bound is that of \cite{R} and is $n^2$; the 
lower bound is $n(n+1)/2$ and is conjectured to be equality.  For $n=2,3,4$ this conjecture, known as Kuzmin's conjecture, has been verified (see \cite{DF}).  In the proof of the 
Nagata-Higman theorem (see \cite{DF}), the bound is computed to be $2^n-1$, which is how $|\mathtt{w}|\leq 7$ in \eqref{generators} arises. 

Let $\xb^*_k$ be the transpose of the matrix of cofactors of $\xb_k$.  In other words, the $(i,j)^{\text{th}}$ entry of
$\xb^*_k$ is $$(-1)^{i+j}\mathrm{Cof}_{ji}(\xb_k);$$ that is, the determinant obtained by
removing the $j^{\text{th}}$ row and $i^{\text{th}}$ column of $\xb_k$.   Let $\M_r^*$ be the monoid generated by $\{\xb_1,\xb_2,...,\xb_r\}$ and
$\{\xb^*_1,\xb^*_2,...,\xb^*_r\}.$  Observe that $$(\xb\yb)^*=\yb^*\xb^*$$ for all $\xb,\yb \in \M^+_r$, and $$\xb\xb^*=\det(\xb)\id.$$  Now let 
$\mathbf{N}_r$ be the normal sub-monoid generated by $$\{\det(\xb_k)\id\ |\ 1\leq k\leq r\},$$ and 
subsequently define $\M_r=\M_r^*/\mathbf{N}_r$.  Notice in $\M_r$, $\xb^*=\xb^{-1}$, and thus $\M_r$ is a group.

We will need the structure of an algebra, and to that end let $\C \M_r$ be the group algebra defined over $\C$ with respect to matrix addition 
and scalar multiplication in $\M_r$.  Likewise, let $\C\M^*_r$ be the semi-group algebra of the monoid $\M^*_r$.

The following commutative diagram relates these objects:
\begin{displaymath}
\begin{CD}
\F_r^+ @>>>\F_r   @= \F_r\\
@VVV  @.         @VVV\\
\M_r^+ @>>> \M_r^* @>>> \M_r\\
@VVV    @VVV       @VVV\\
\C \M_r^+ @>>>\C \M_r^* @>>> \C \M_r @>\mathrm{tr}>> \C[\X].\\
\end{CD} 
\end{displaymath}

\subsection{Relations}

The Cayley-Hamilton theorem applies to this context and so for any $\xb\in \C \M_r$, 
\begin{align} \label{cayham}
\xb^3-\tr{\xb}\xb^2+ \tr{ \xb^* } \xb-\det(\xb)\id=0.
\end{align}

By direct calculation, or by Newton's trace formulas

\begin{align}
\tr{ \mathbf{x}^* }&=\frac{1}{2}\left(\tr{ \mathbf{x} }^2-\tr{ \mathbf{x}^2 }\right)\label{trinv}.
\end{align}

Together \eqref{cayham} and \eqref{trinv} imply
\begin{align}
\det(\mathbf{x})&=\frac{1}{3}\tr{\xb^3}+\frac{1}{6}\tr{\xb}^3-\frac{1}{2}\tr{\xb}\tr{\xb^2}.\label{dettr}
\end{align}

\begin{rem}
In general the coefficients of the characteristic polynomial of an $n\times n$ matrix are the elementary symmetric polynomials in the eigenvalues of the matrix.  By Newton's
formulas these are trace expressions in powers of the matrix.  So one may use this and the general method of polarization, which we demonstrate below, to develop trace identities for larger 
size matrices.
\end{rem}

Computations similar to those that follow may be found in \cite{MS, SR}; the process is standard and is generally known as (partial) polarization, or multilinearization.  
For any $\xb, \yb 
\in \C \M_r$ and any $\lambda \in \C$, equation 
$\eqref{cayham}$ implies
\begin{align}\label{cayhamsum}
(\xb+\lambda \yb)^3-\tr{\xb+\lambda \yb}(\xb+\lambda \yb)^2 + \tr{(\xb+\lambda 
\yb)^*}(\xb+\lambda\yb)-\det(\xb+\lambda \yb)\id=0.
\end{align}

Using equations $\eqref{cayham}$, $\eqref{trinv}$, and $\eqref{dettr}$, equation $\eqref{cayhamsum}$ simplifies to
\begin{align*}
0=&\lambda^2\bigg(\xb\yb^2+\yb^2\xb+\yb\xb\yb-\tr{\xb}\yb^2-\tr{\yb}\xb\yb-\tr{\yb}\yb\xb 
+\frac{1}{2}\tr{\yb}^2\xb-\frac{1}{2}\tr{\yb^2}\xb+\\
&\tr{\xb}\tr{\yb}\yb-\tr{\xb\yb}\yb -\tr{\xb\yb^2}\id-\frac{1}{2}\tr{\xb}\tr{\yb}^2\id
+\frac{1}{2}\tr{\xb}\tr{\yb^2}\id +\tr{\yb}\tr{\xb\yb}\id\bigg)+\\
&\lambda\bigg(\yb\xb^2+\xb^2\yb+\xb\yb\xb-\tr{\yb}\xb^2-\tr{\xb}\yb\xb-\tr{\xb}\xb\yb
+\frac{1}{2}\tr{\xb}^2\yb-\frac{1}{2}\tr{\xb^2}\yb+\\
&\tr{\xb}\tr{\yb}\xb-\tr{\xb\yb}\xb
-\tr{\yb\xb^2}\id-\frac{1}{2}\tr{\yb}\tr{\xb}^2\id+\frac{1}{2}\tr{\yb}\tr{\xb^2}\id 
+\tr{\xb}\tr{\xb\yb}\id\bigg). 
\end{align*}

Thus, by Vandermonde arguments (see \cite{Ro}) we have the partial polarization of $\eqref{cayham}$
\begin{align}\label{polarization}
&\yb\xb^2 +\xb^2\yb+\xb\yb\xb= \tr{\yb}\xb^2+\tr{\xb}\yb\xb +\tr{\xb}\xb\yb 
-\tr{\xb}\tr{\yb}\xb +\tr{\xb\yb}\xb+  \nonumber\\ 
&\tr{\yb\xb^2}\id-\tr{\xb}\tr{\xb\yb}\id -\frac{1}{2}\left(\tr{\xb}^2\yb -\tr{\xb^2}\yb- \tr{\yb}\tr{\xb}^2\id +\tr{\yb}\tr{\xb^2}\id\right).
\end{align}

Define $\pol(\xb,\yb)$ to be the right hand side of equation $\eqref{polarization}$;  that is, 
\begin{align}\pol(\xb,\yb)=\yb\xb^2 +\xb^2\yb+\xb\yb\xb.\label{polarization:2}\end{align} 

Then substituting $\xb$ by the sum $\xb + \zb$ in equation $\eqref{polarization:2}$, yields the full polarization of $\eqref{cayham}$ 
\begin{align}\label{fundamental}
\xb\zb\yb+\zb\xb\yb+\yb\xb\zb+\yb\zb\xb+\xb\yb\zb+\zb\yb\xb=\pol(\xb+\zb,\yb)-\pol(\xb,\yb)-\pol(\zb,\yb).
\end{align}

If $\xb, \yb \in \M_r$ then multiplying equation $\eqref{cayham}$ on the right by $\xb^{-1}\yb$ yields,
\begin{eqnarray}\label{cayham2}
\xb^2\yb-\tr{\xb}\xb \yb +\tr{\xb^{-1}}\yb-\xb^{-1}\yb=0.
\end{eqnarray}

Suppose $\xb,\yb\in \M_r$.  Multiplying equation $\eqref{polarization}$ on the left by $\yb^{-1}\xb^{-1}$ and on the right by $\xb^{-1}$, followed by 
taking the trace, and using equation $\eqref{cayham2}$ appropriately, provides the commutator trace relation 
\begin{align}\label{polyp1}
\tr{\xb\yb\xb^{-1}\yb^{-1}}=& -\tr{\yb\xb\yb^{-1}\xb^{-1}}+\tr{\xb}\tr{\xb^{-1}}\tr{\yb}\tr{\yb^{-1}} \nonumber\\
                            &+\tr{\xb}\tr{\xb^{-1}} +\tr{\yb}\tr{\yb^{-1}}+\tr{\xb\yb}\tr{\xb^{-1}\yb^{-1}}\nonumber\\
                            &+\tr{\xb\yb^{-1}}\tr{\xb^{-1}\yb}-\tr{\xb^{-1}}\tr{\yb}\tr{\xb\yb^{-1}} \\
                            &-\tr{\xb}\tr{\yb^{-1}}\tr{\xb^{-1}\yb} -\tr{\xb}\tr{\yb}\tr{\xb^{-1}\yb^{-1}}\nonumber\\ 
                            &-\tr{\xb\yb}\tr{\xb^{-1}}\tr{\yb^{-1}}-3.\nonumber
\end{align}

\subsection{Generators}
From $\eqref{generators}$, we need only consider words in $\F^+_r$ of length $7$ or less.  In \cite{SR} it is shown that this length may be taken to be $6$.  We give a similar 
argument here since the development of the result provides many useful relations, and a constructive algorithm for word reduction that is of computational significance.  

The length of a reduced word is defined to be the number of letters, counting multiplicity, in the word.  We now define 
the {\it weighted length}, denoted by $|\cdot|_w$, to be the number of letters of a reduced word having 
positive exponent plus twice the number of letters having negative exponent, again counting multiplicity.

For example, in $\F_2$, we have $|\xt_1\xt_2|=|\xt_1\xt_2|_w=2$ but $|\xt_1^3\xt_2^{-2}|=3+2=5$ while
$|\xt_1^3\xt_2^{-2}|_w=3+2\cdot 2=7$.

For a polynomial expression $e$ in generic matrices with coefficients in $\C[\X]$, we define the {\it degree of $e$}, denoted by $|\!|e|\!|$, to be the largest  
weighted length of monomial words in the expression of $e$ that is minimal among all such expressions for $e$.  Additionally, we define the {\it trace degree of 
$e$}, denoted by $|\!|e|\!|_{\mathrm{tr}}$, to be the maximal degree over all monomial words within a trace coefficient of $e$.

For example, when $\xb, \yb\in \M_r$, $|\!|\pol(\xb,\yb)|\!|\leq \mathrm{max}\{2|\!|\xb|\!|,|\!|\xb|\!|+|\!|\yb|\!|\},$ while 
$|\!|\pol(\xb,\yb)|\!|_{\mathrm{tr}}\leq 2|\!|\xb|\!|+|\!|\yb|\!|$.

We remark that given two such expressions $e_1$ and $e_2$, $$|\!|e_1e_2|\!|\leq|\!|e_1|\!|+|\!|e_2|\!|\ \text{ and }\ 
|\!|e_1e_2|\!|_{\mathrm{tr}}\leq\mathrm{max}\{|\!|e_1|\!|_{\mathrm{tr}},|\!|e_2|\!|_{\mathrm{tr}}\}.$$  

We are now prepared to characterize the generators of $\C[\X]$.

\begin{lem}
$\C[\R]^\G$ is generated by $\tr{\mathbf{w}}$ such that $\mathtt{w}\in\F_r$ is cyclicly reduced, $|\mathtt{w}|_w \leq 6$, and all 
exponents of letters in $\mathtt{w}$ are $\pm 1$.
\end{lem}

\begin{proof}
For $n\geq 2$, equations $\eqref{cayham}$ and $\eqref{cayham2}$ determine equation 
\begin{align}\label{powerreduce}
\tr{\ub\xb^n\vb}=&\tr{\xb}\tr{\ub\xb^{n-1}\vb}-\tr{\xb^{-1}}\tr{\ub\xb^{n-2}\vb}+\tr{\ub\xb^{n-3}\vb},
\end{align}
which recursively reduces $\tr{\mathbf{w}}$ to a polynomial in traces of words having no letter with exponent 
other than $\pm 1$.  If however $n\leq -2$ then we first apply equation \eqref{cayham2} and then use \eqref{powerreduce}.  Hence it follows 
that $\mathtt{w}$ may be taken to be cyclically reduced, having all letters with exponent $\pm 1$.

It remains to show that the word length may be taken to be less than or equal to $6$.

Substituting $\xb\mapsto \yb$ and $\yb\mapsto \xb\zb$ in equation $\eqref{polarization:2}$, and multiplying the resulting expression on the left by $\xb$ gives
\begin{align}\label{eq:4}\xb^2\zb\yb^2=-(\xb\yb^2\xb)\zb-(\xb\yb\xb)\zb\yb+\xb\pol(\yb,\xb).\end{align}  

Replacing $\yb\mapsto \yb^2$ in equation $\eqref{polarization:2}$ produces $$\yb^2\xb^2+\xb^2\yb^2+\xb\yb^2\xb=\pol(\xb,\yb^2),$$ which substituted into equation $\eqref{eq:4}$ 
yields equation \begin{align}\label{eq:5}\xb^2\zb\yb^2=(\yb^2\xb^2+\xb^2\yb^2-\pol(\xb,\yb^2))\zb+(\yb\xb^2+\xb^2\yb-\pol(\xb,\yb))\zb\yb+\xb\pol(\yb,\xb\zb).\end{align}

Now substituting $$\pol(\yb,\xb^2\zb)=\xb^2\zb\yb^2+\yb^2\xb^2\zb+\yb\xb^2\zb\yb,$$ and $$\xb^2\pol(\yb,\zb)=\xb^2\zb\yb^2+\xb^2\yb^2\zb+\xb^2\yb\zb\yb$$
into equation $\eqref{eq:5}$ results in 
\begin{align}3\xb^2\zb\yb^2=\pol(\yb,\xb^2\zb)+\xb\pol(\yb,\xb\zb)-\pol(\xb,\yb^2)\zb-\pol(\xb,\yb)\zb\yb+\xb^2\pol(\yb,\zb)\label{eq:6}.\end{align}

Thus, $$|\!|\xb^2\zb\yb^2|\!|< 2|\!|\xb|\!|+|\!|\zb|\!|+2|\!|\yb|\!|\ \text{ and } \tdg{\xb^2\zb\yb^2}\leq 2|\!|\xb|\!|+|\!|\zb|\!|+2|\!|\yb|\!|.$$

\begin{rem}
In the proof of the Nagata-Higman theorem, the two-sided ideal of polynomial trace relations, for $n=3$, is shown to be generated as a vector space by $\pol(\ub,\vb)$, $\ub^3$, and 
equation \eqref{fundamental} evaluated at monomial words $\ub$, $\vb$, and $\wb$.  Equation $\eqref{eq:6}$ shows $\xb^2\zb\yb^2$ is in this ideal, and consequently its degree is less 
than its word length.  However, one can conclude $\xb^2\zb\yb^2$ is in this ideal from more general considerations and avoid the above calculation $($see \cite{DF}, page $76)$.  
\end{rem}

For the remainder of the argument assume $\xb,\yb,\zb,\ub,\vb,\wb$ are of length $1$.  Replacing $\yb\mapsto \ub+\vb$ in equation $\eqref{eq:6}$ we deduce 
$|\!|\xb^2\zb(\ub^2+\ub\vb+\vb\ub+\vb^2)|\!|\leq 4.$  This in turn implies $|\!|\xb^2\zb(\ub\vb+\vb\ub)|\!|\leq 4$ and so $|\!|\xb^2\zb\wb(\ub\vb+\vb\ub)|\!|\leq 5.$  In a 
like manner, we have that both $\dg{\xb^2\zb(\wb\ub\vb+\vb\wb\ub)}\leq 5$ and $\dg{\xb^2\zb(\wb\vb+\vb\wb)\ub}\leq 5$.  Hence we conclude that 
$$\dg{2\xb^2\zb\wb\ub\vb}=\dg{\xb^2\zb\wb(\ub\vb+\vb\ub)+\xb^2\zb(\wb\ub\vb+\vb\wb\ub)-\xb^2\zb(\wb\vb+\vb\wb)\ub}\leq 5,$$ and $$\tdg{2\xb^2\zb\wb\ub\vb}\leq 6.$$

Replacing $\xb\mapsto \xb +\yb$ in $\xb^2\zb\wb\ub\vb$ we come to the conclusion that $\dg{\xb\yb\zb\wb\ub\vb+\yb\xb\zb\wb\ub\vb}\leq 5$.  In other words, permuting $\xb$ and $\yb$  
introduces a factor of $-1$ and a polynomial term of lesser degree.  Slight variation in our analysis concludes the same result for any transposition of two adjacent letters in the word
$\xb\yb\zb\wb\ub\vb$.

Therefore, if $\sigma$ is a permutation of the letters $\xb,\yb,\zb,\ub,\vb,\wb$ then 
$$\dg{\xb\yb\zb\ub\vb\wb+\mathrm{sgn}(\sigma)\sigma(\xb\yb\zb\ub\vb\wb)}\leq 5\ \text{while}\ \tdg{\xb\yb\zb\ub\vb\wb+\mathrm{sgn}(\sigma)\sigma(\xb\yb\zb\ub\vb\wb)}\leq 6.$$

Lastly, making the substitutions $\xb\mapsto \xb\yb$, $\yb\mapsto \zb\ub$, and $\zb\mapsto \vb\wb$ in the fundamental expression $\eqref{fundamental}$, we derive
\begin{align} \label{eq:7} {\xb\yb}\vb\wb{\zb\ub}+&\vb\wb{\xb\yb}{\zb\ub}+{\zb\ub}{\xb\yb}\vb\wb+{\zb\ub}\vb\wb{\xb\yb} 
+{\xb\yb}{\zb\ub}\vb\wb+\vb\wb{\zb\ub}{\xb\yb}=\nonumber\\&\pol({\xb\yb}+\vb\wb,{\zb\ub})-\pol({\xb\yb},\zb\ub)-\pol(\vb\wb, {\zb\ub}).\end{align}

However, each word on the left hand side of equation $\eqref{eq:7}$ is an even permutation of the first, so 

$$\dg{6{\xb\yb}\vb\wb{\zb\ub}}\leq 5\ \text{and}\ \tdg{6{\xb\yb}\vb\wb{\zb\ub}}\leq 6.$$  
Hence, if $\wb$ is a word of length $7$ or more, then $\tdg{\tr{\wb}}\leq 6$.  
Moreover, this process gives an iterative algorithm for reducing such an expression.
\end{proof}

As an immediate consequence we have the following description of sufficient generators of $\C[\X]$.  

\begin{cor}\label{generatorform}
$\C[\X]$ is generated by traces of the form 
\begin{align*}&\tr{\xb_i},\tr{\xb_i^{-1}}, \tr{\xb_i\xb_j}, \tr{\xb_i\xb_j\xb_k}, \tr{\xb_i\xb_j^{-1}},\tr{\xb_i^{-1}\xb_j^{-1}},\tr{\xb_i\xb_j\xb_k^{-1}},\\
&\tr{\xb_i\xb_j\xb_k\xb_l},\tr{\xb_i\xb_j\xb_k\xb_l\xb_m},\tr{\xb_i\xb_j\xb_k\xb_l^{-1}}, \tr{\xb_i\xb_j\xb_k\xb_j^{-1}},\tr{\xb_i\xb_j^{-1}\xb_k^{-1}},\\
& \tr{\xb_i^{-1}\xb_j^{-1}\xb_k^{-1}}, \tr{\xb_i\xb_j\xb_k^{-1}\xb_l^{-1}},\tr{\xb_i\xb_j\xb_k^{-1}\xb_j^{-1}}, 
\tr{\xb_i\xb_j\xb_i^{-1}\xb_j^{-1}},\\
&\tr{\xb_i\xb_j\xb_k\xb_l\xb_m^{-1}}, \tr{\xb_i\xb_j\xb_k\xb_l\xb_k^{-1}}, \tr{\xb_i\xb_j\xb_k\xb_l\xb_j^{-1}}, \tr{\xb_i\xb_j\xb_k\xb_l\xb_m\xb_n},\end{align*} where 
$1\leq i\not=j\not=k\not=l\not=m\not=n \leq r$.  
\end{cor}

\begin{proof}
First, consider generators of type $\tr{\ub^{-1}\wb\xb^{-1}\zb}$.  It can be 
shown that $\tr{\ub\vb\wb\xb\yb\zb}+\tr{\ub\vb\wb\yb\xb\zb}+\tr{\vb\ub\wb\xb\yb\zb}+\tr{\vb\ub\wb\yb\xb\zb}$ has trace degree $5$.  Setting $\ub=\vb$ and $\xb=\yb$ and subsequently 
interchanging words with squares to those with inverses, we find generators of the form $\tr{\ub^{-1}\wb\xb^{-1}\zb}$ can be freely eliminated; that is, inverses can be assumed to be adjacent.

It remains to show that letters may be taken to be distinct.  Equation $\eqref{polarization:2}$ implies that for any letter $\xb$ and any monomial words $\wb_1, 
\wb_2, \wb_3$,
$$\tr{\wb_1\xb^{\pm 1}\wb_2 \xb^{\pm 1}\wb_3}=-\tr{\wb_1\xb^{\pm 2}\wb_2 \wb_3}-\tr{\wb_1\wb_2 \xb^{\pm 2}\wb_3}+\tr{\wb_1\pol(\xb^{\pm 
1}, \wb_2)\wb_3}.$$  However, by subsequently reducing the words having letters with exponent not $\pm 1$, we conclude that expressions of the form $\tr{\wb_1\xb^{\pm 
1}\wb_2 \xb^{\pm 1}\wb_3}$ are unnecessary.
\end{proof}

This result can be refined using the work of \cite{AP}, where explicit {\it minimal} generators are formulated in a more general context.  In an upcoming paper, we will address 
the issue of minimality for our generators, as well as provide a maximal subset that is algebraically independent.  This subset will allow for a generalization of the 
symmetry described in section $\ref{outer}$.

\section{Structure of $\C[\G\times \G]^{\G}$}
\subsection{Minimal Generators}
As a consequence of Corollary $\ref{generatorform}$, we have
\begin{lem}\label{gens}
$\C [\G\times \G]^\G$ is generated by 
\begin{align*}
&\tr{\xb_1},\ \ \tr{\xb_2},\ \ \tr{\xb_1\xb_2},\ \ \tr{\xb_1\xb_2^{-1}},\ \ \tr{\xb_1^{-1}},\\
&\tr{\xb_2^{-1}},\ \ \tr{\xb_1^{-1}\xb_2^{-1}},\ \ \tr{\xb_1^{-1}\xb_2},\ \ \tr{\xb_1\xb_2\xb_1^{-1}\xb_2^{-1}}.
\end{align*}
\end{lem}

\begin{proof}
The words of weighted length $1,2,3,4$ with exponents $\pm 1$ are unambiguously cyclically equivalent to one of $$\xt_1, 
\xt_2, \xt_1\xt_2, \xt_1\xt_2^{-1}, \xt_2\xt_1^{-1}, \xt_1^{-1}\xt_2^{-1}, (\xt_1\xt_2)^2.$$  But equation $\eqref{cayham2}$ reduces the 
latter most of these in terms of the others.  All words in two letters of length $5$ are cyclically equivalent to a 
word with an exponent whose magnitude is greater than $1$, except $\xt_1\xt_2^{-1}\xt_1\xt_2$, and $\xt_2\xt_1^{-1}\xt_2\xt_1$.  
Both are cyclically equivalent to $(\xt_i\xt_j)^2\xt_j^{-2}$ which in turn, by equation $\eqref{powerreduce}$ reduces to 
expressions in the other variables.  The only words of weighted length $6$ and with 
exponents only $\pm 1$ are $\xt_1\xt_2 \xt_1^{-1}\xt_2^{-1}$, its inverse, and $(\xt_1\xt_2)^3$.  But the latter most of 
these is reduced by equation $\eqref{cayham}$.  Lastly, letting $\xb=\xb_1$ and $\yb=\xb_2$ in equation 
$\eqref{polyp1}$, we have
\begin{align}\label{polyp2}
\tr{\xb_2\xb_1\xb_2^{-1}\xb_1^{-1}}=& -\tr{\xb_1\xb_2\xb_1^{-1}\xb_2^{-1}}+\tr{\xb_1}\tr{\xb_1^{-1}}\tr{\xb_2}\tr{\xb_2^{-1}}\nonumber\\ 
                                    &+\tr{\xb_1}\tr{\xb_1^{-1}}+\tr{\xb_2}\tr{\xb_2^{-1}}+\tr{\xb_1\xb_2}\tr{\xb_1^{-1}\xb_2^{-1}}\nonumber\\
                                    &+\tr{\xb_1\xb_2^{-1}}\tr{\xb_1^{-1}\xb_2}-\tr{\xb_1^{-1}}\tr{\xb_2}\tr{\xb_1\xb_2^{-1}}\\
                                    &-\tr{\xb_1}\tr{\xb_2^{-1}}\tr{\xb_1^{-1}\xb_2} -\tr{\xb_1}\tr{\xb_2}\tr{\xb_1^{-1}\xb_2^{-1}}\nonumber\\
                                    &-\tr{\xb_1\xb_2}\tr{\xb_1^{-1}}\tr{\xb_2^{-1}}-3, \nonumber
\end{align}
which expresses the trace of the inverse of the commutator in terms of the other expressions.
\end{proof}

\subsection{$\mathbb{Z}_3^{\times 2}$-Grading}
The center of $\G$ is $\zeta(\G)=\{\omega\id\ |\ \omega^3=1\}\cong\mathbb{Z}_3$.  There is an action of $\zeta(\G)^{\times 2}$ on $\C[\X]$ given by
$$(\omega_1\id,\omega_2\id)\cdot\tr{\mathbf{w}(\xb_1,\xb_2)}=\tr{\mathbf{w}(\omega_1\xb_1,\omega_2\xb_2)}=
\omega_1^{|\mathbf{w}(\xb_1,\id)|_w}\omega_2^{|\mathbf{w}(\id,\xb_2)|_w}\tr{\mathbf{w}(\xb_1,\xb_2)}.$$  Applying this action to the generators and recording the
orbit by a $9$-tuple, all generators are distinguished.  Consequently, we have

\begin{prop}
$$\C[\X]=\sum_{(\omega_1,\omega_2)\in\mathbb{Z}_3\times \mathbb{Z}_3} \C[\X]_{(\omega_1,\omega_2)}$$ is a $\mathbb{Z}_3\times \mathbb{Z}_3$-graded ring.  The summand 
$\C[\X]_{(\omega_1,\omega_2)}$ is the linear span over $\C$ of all monomials whose orbit under $\mathbb{Z}_3\times \mathbb{Z}_3$ equals one of the orbits of the 
nine orbit types corresponding to the minimal generators. 
\end{prop}

In fact the situation is general.  For a rank $r$ free group, $\mathbb{Z}_3^{\times r}$ acts on the generators of $\C[\X]$ and gives a filtration.  However,
since the relations are polarizations of the Cayley-Hamilton polynomial, which itself has a zero grading, no relation can compromise summands.  So the filtration is a
grading.

\subsection{Hypersurface in $\C^9$}

Let $$\overline{R}=\C[t_{(1)},t_{(-1)},t_{(2)},t_{(-2)},t_{(3)},t_{(-3)},t_{(4)},t_{(-4)}, 
t_{(5)},t_{(-5)}]$$ be the complex polynomial ring freely generated by $\{t_{(\pm i)},\ 1\leq i\leq 5\},$ and let 
$$R=\C[t_{(1)},t_{(-1)},t_{(2)},t_{(-2)},t_{(3)},t_{(-3)},t_{(4)},t_{(-4)}]$$ be its subring generated 
by $\{t_{(\pm i)},\  1\leq i\leq 4\},$  so 
$\overline{R}=R[t_{(5)},t_{(-5)}].$
Define the following ring homomorphism,
$$R[t_{(5)},t_{(-5)}]\stackrel{\tiny \Pi}{\longrightarrow} \C[\G\times \G]^\G$$ by

\begin{center}
\begin{tabular}{ll}
$t_{(1)}\mapsto\tr{\mathbf{x}_1}$ & $t_{(-1)}\mapsto\tr{\mathbf{x}_1^{-1}}$\\
$t_{(2)}\mapsto\tr{\mathbf{x}_2}$& $t_{(-2)}\mapsto\tr{\mathbf{x}_2^{-1}}$\\
$t_{(3)}\mapsto\tr{\mathbf{x}_1\mathbf{x}_2}$& $t_{(-3)}\mapsto\tr{\mathbf{x}_1^{-1}\mathbf{x}_2^{-1}}$\\
$t_{(4)}\mapsto\tr{\mathbf{x}_1\mathbf{x}_2^{-1}}$& $t_{(-4)}\mapsto\tr{\mathbf{x}_1^{-1}\mathbf{x}_2}$\\
$t_{(5)}\mapsto\tr{\mathbf{x}_1\mathbf{x}_2\mathbf{x}_1^{-1}\mathbf{x}_2^{-1}}$& $t_{(-5)}\mapsto\tr{\mathbf{x}_2\mathbf{x}_1\mathbf{x}_2^{-1}\mathbf{x}_1^{-1}}$.
\end{tabular}
\end{center}

It follows from Lemma \ref{gens} that $$\C[\X]\cong R[t_{(5)}, t_{(-5)}]/\ker(\Pi).$$ In other words, $\Pi$ is a 
surjective algebra morphism.\\
We define $$P=t_{(1)}t_{(-1)}t_{(2)}t_{(-2)}-t_{(1)}t_{(2)}t_{(-3)}-t_{(-1)}t_{(-2)}t_{(3)}
-t_{(1)}t_{(-2)}t_{(-4)}-t_{(-1)}t_{(2)}t_{(4)}$$ 
$$+t_{(1)}t_{(-1)}+t_{(2)}t_{(-2)}+t_{(3)}t_{(-3)}+t_{(4)}t_{(-4)}-3,$$ and so $P\in R$.  
Moreover, by equation $\eqref{polyp2}$, $$P-(t_{(5)}+t_{(-5)}) \in \ker(\Pi).$$

Hence it follows that the composite map $$R[t_{(5)}]\hookrightarrow R[t_{(5)}, t_{(-5)}]\twoheadrightarrow R[t_{(5)}, 
t_{(-5)}]/\ker(\Pi),$$ is an epimorphism.  Let $I$ be the kernel of this composite map, and suppose there exists $Q\in 
R$ so $Q-t_{(5)}t_{(-5)}\in \ker(\Pi)$ as well.

Then under this assumption, we prove
\begin{lem}\label{ideallemma}
$I$ is principally generated by the polynomial 
\begin{equation}\label{ideal}
t_{(5)}^2-Pt_{(5)}+Q.
\end{equation}
\end{lem}
\begin{proof}
The following argument is an adaptation of one found in \cite{N}.

Certainly, $t_{(5)}^2-Pt_{(5)}+Q\in I$ for it maps into $R[t_{(5)},t_{(-5)}]/\ker(\Pi)$ to the coset 
representative $t_{(5)}^2-(t_{(5)}+t_{(-5)})t_{(5)}+t_{(5)}t_{(-5)}=0$.

On the other hand, observe $$R[t_{(5)}]/I\cong R[t_{(5)}, t_{(-5)}]/\ker(\Pi)\cong \C[\X],$$ the
dimension of $\X$ is $8$, and $R[t_{(5)}]/I$ has at most $9$ generators.   Then it must be the case that $I$ is principally generated since $R[\ti{5}]$ is a
U.F.D., and thus a co-dimension $1$ irreducible subvariety of $\C^9$ must be given by one equation (see \cite{S} page $69$).  Moreover, $I$ is
non-zero since otherwise the resulting dimension would necessarily be too large.

Seeking a contradiction, suppose there exists a polynomial identity comprised of only elements of $R$.  Then Krull's dimension theorem (see page $68$ in \cite{S}) implies 
$t_{(5)}$ is free.  In other words, given any restriction of the generators of $R$, $t_{(5)}$ is not determined.  Consider
$\mathrm{SL}_3(\mathrm{SL}(2,\C))^{\times 2}\subset \G^{\times 2}$; that is, matrices of the form $\left(
\begin{array}{ccc}
a & b & 0\\
c & d & 0\\
0& 0 & 1\\
\end{array}\right)$ so $ad-bc=1$.  Then by restricting to pairs of such matrices, we deduce that $$\tr{\xb_1\xb_2\xb_1^{-1}\xb_2^{-1}}=
\tr{\xb_2\xb_1\xb_2^{-1}\xb_1^{-1}},$$ since for all $\xb\in \mathrm{SL}(2,\C),$ $\tr{\xb}=\tr{\xb^{-1}}$.
Then equation $\eqref{polyp2}$ becomes $$t_{(5)}=P/2,$$ which is decidedly not free of the 
generators of $R$.  Thus, the generators of $R$ are algebraically independent in $R[t_{(5)}]/I$.  

Since $I$ is principal and contains a monic quadratic over $R$, its generator is expression $\eqref{ideal}$, or a 
factor thereof.  We have just showed that there are no degree zero relations, with respect to $t_{(5)}$.  However, 
if $I$ is generated by a linear polynomial over $R$ then $t_{(5)}$ is determined by the generators of 
$R$ alone.  However this in turn would imply that all representations who agree by evaluation in $R$ also agree by evaluation under $t_{(5)}$.

Consider the representations 

\begin{tabular}{ccc}
$\F_2 \stackrel{\rho_1}{\longrightarrow} \G$ & &$\F_2 \stackrel{\rho_2}{\longrightarrow} \G$\\
$\xt_1 \longmapsto
\left(
\begin{array}{ccc}
a & 0 & 0\\
0 & b & 0\\
0& 0 & 1/ab\\
\end{array}\right)$ &\ \text{and}\  &
$\xt_1 \longmapsto
\left(
\begin{array}{ccc}
a & 0 & 0\\
0 & b & 0\\
0& 0 & 1/ab\\
\end{array}\right)$\\
$\xt_2 \longmapsto
\frac{1}{4^{1/3}}\left(
\begin{array}{ccc}
1 & 1 & -1\\
1 & -1 & 1\\
-1& -1 & -1\\
\end{array}\right)$
& & $\xt_2 \longmapsto
\frac{1}{4^{1/3}}\left(
\begin{array}{ccc}
1 &-1 & 1\\
-1 & -1 & -1\\
1& 1 & -1\\
\end{array}\right).$
\end{tabular}

It is a direct calculation to verify that they agree upon evaluation in $R$ but disagree under $t_{(5)}$.
\end{proof}

Lemmas \ref{gens} and \ref{ideallemma} together imply the following theorem whose result, in part, was given by \cite{We}, and later by \cite{Si}, and may also be inferred by the work of 
\cite{T}.

\begin{thm}\label{ranktwo}
$\G^{\times 2}\aq\G$ is isomorphic to a degree $6$ affine hyper-surface in $\C^9$, which maps onto $\C^8$.
\end{thm}

\begin{proof}
The degree of $Q$ will be apparent when we explicitly write it down.  It remains to show that $\X\rightarrow 
\C^8$ is a surjection.  To this end, let 
$(z_1-\zeta_1,...,z_8-\zeta_8)$ be a maximal ideal in the coordinate ring of $\C^8$.  
Moreover, let $\zeta_9$ be defined to be a solution to $t^2-P(\zeta_1,...,\zeta_8)t+Q(\zeta_1,...,\zeta_8)=0$.  Then 
$(t_{(1)}-\zeta_1,t_{(-1)}-\zeta_2,...,t_{(-4)}-\zeta_8,t_{(5)}-\zeta_9)+I$ is a maximal 
ideal in $\C[\X]$, and so all maximal ideals of $\C[\C^8]$ are images of such in $\C[\X]$. 
\end{proof}

\subsection{Singular Locus of $\X$.}

The surjection $\X\to \C^8$ is generically $2$-to-$1$; that is, there are exactly two solutions to $$t^2-P(\zeta_1,...,\zeta_8)t+Q(\zeta_1,...,\zeta_8)=0$$ for every
point in $\C^8$ except where $P^2-4Q=0$.  In this case, $$0=(t_{(5)}+t_{(-5)})^2-4t_{(5)}t_{(-5)}=(t_{(5)}-t_{(-5)})^2$$ which implies $t_{(5)}=t_{(-5)}=P/2$.  On the 
other hand, at points in $\X$ if $t_{(5)}=P/2$, then $P^2-4Q=0$.  Let $\frak{L}$ denote the locus of solutions to $P^2-4Q=0$ in $\X$, which is a closed subset 
of $\X$.

It is readily observed that the partial derivative with respect to $\ti{5}$ of $t_{(5)}^2-Pt_{(5)}+Q$ is zero if and only if $t_{(5)}=P/2$.  The singular set in $\X$, 
denoted by $\frak{J}$, is the closed subset cut out by the Jacobian ideal; that is, the ideal generated by the formal partial derivatives of $\ti{5}^2-P\ti{5}+Q$.  
Thus $\frak{J}\subset\frak{L}$.  

If $\frak{H}\hookrightarrow \G$ is a sub-algebraic group, then we define $\frak{H}^{\times r}\aq \G$ to be the image of $\frak{H}^{\times r}\hookrightarrow \R \to \X$.

In the proof of Lemma $\ref{ideallemma}$, we observed $\mathrm{SL}_3(\mathrm{SL}(2,\C))^{\times 2}\aq \G \subset
\frak{L}$.  Additionally, since matrices of the form
$\left(
\begin{array}{ccc}
a & 0 & 0\\
0 & b & 0\\
0& 0 & 1/ab\\
\end{array}\right)$ commute, restricting to pairs of such matrices enforces the relation
$$\tr{\xb_1\xb_2\xb^{-1}_1\xb^{-1}_2}=3=\tr{\xb_2\xb_1\xb^{-1}_2\xb^{-1}_1}.$$  Let $\mathrm{SL}_3(\C^*\times \C^*)$ denote the subset of such matrices in $\G$.
Consequently, $\mathrm{SL}_3(\C^*\times \C^*)^{\times 2}\aq \G\subset \frak{L}$ as well.  We claim both sets satisfy all the generators of the Jacobian ideal, and so 
are singular.  The Jacobian ideal is generated by the polynomials $-\ti{5}\dP{i}+\dQ{i}$ for $1\leq |i|\leq 4$, and $2\ti{5}-P$.  
Using the formulas for $P$ and $Q$ (see section \ref{detq}), we explicitly write out the generators of the Jacobian 
ideal (see \cite{L} for details).  Then evaluating these polynomials 
at pairs of generic matrices in either $\mathrm{SL}_3(\mathrm{SL}(2,\C))$ or $\mathrm{SL}_3(\C^*\times \C^*)$ we verify 
that all partials vanish using {\it Mathematica} (\cite{W}).  It turns out these examples are prototypical.

Let $\mathrm{SL}_3(\mathrm{GL}(2,\C))$ be the subset of $\G$ consisting of elements of the form 
$\left(
\begin{array}{ccc}
a & b & 0\\
c & d & 0\\
0& 0 & \frac{1}{ad-bc}
\end{array}\right).$

Notice that $\mathrm{SL}_3(\C^*\times \C^*)^{\times 2}\aq \G$ and $\mathrm{SL}_3(\mathrm{SL}(2,\C))^{\times 2}\aq \G$ are contained in 
$\mathrm{SL}_3(\mathrm{GL}(2,\C))^{\times 2}\aq \G$.  Again, using {\it Mathematica} we evaluate all generating polynomials of the Jacobian ideal on pairs of generic 
matrices in $\mathrm{SL}_3(\mathrm{GL}(2,\C))$.  Since all partials vanish, $\mathrm{SL}_3(\mathrm{GL}(2,\C))^{\times 2}\aq \G$ is singular in $\X$ as well. 

In general, if $[\rho] \in \G^{\times r}\aq \G$ is singular, then its orbit has positive-dimensional isotropy.  Any {\it completely reducible} representation (these 
parameterize $\G^{\times r}\aq \G$ as an orbit space) that is not {\it irreducible} is conjugate to an element in $\mathrm{SL}_3(\mathrm{GL}(2,\C))^{\times r}$.  
This follows since there must be a shared eigenvector with respect to its generic matrices, if the representation reduces at all.  Irreducible representations are known to be 
non-singular, and their isotropy is zero dimensional.  Consequently, it follows that in general the singular set of $\G^{\times r}\aq \G$ is contained in 
$\mathrm{SL}_3(\mathrm{GL}(2,\C))^{\times r}\aq \G$.  

In the case of a free group of rank $1$, there are no singular points in the quotient and so the identity, which has maximal isotropy, remains non-singular.  Hence the 
converse inclusion does not generally hold.  In the case of a free group of rank $2$, the situation is much better.  In fact, we have already established

\begin{thm}\label{singulartheorem}
A completely reducible representation in $\G^{\times 2}\aq \G$ is singular if and only if its orbit has positive-dimensional isotropy; that is, 
$\frak{J}=\mathrm{SL}_3(\mathrm{GL}(2,\C))^{\times 2}\aq \G$.
\end{thm}  

As a final note, we give an example of a non-singular representation in the branching locus (actually we give a $2$-dimensional family in $\frak{L}-\frak{J}$):

\begin{tabular}{c}
$\F_2 \stackrel{\rho}{\longrightarrow} \G$ \\
$\xt_1 \longmapsto
\left(
\begin{array}{ccc}
a & 0 & 0\\
0 & a & 0\\
0& 0 & 1/a^2\\
\end{array}\right)$

$\xt_2 \longmapsto
\frac{c^{1/3}}{4^{1/3}}\left(
\begin{array}{ccc}
1 & 1 & -1\\
1 & -1 & 1\\
-1/c& -1/c & -1/c\\
\end{array}\right),$
\end{tabular}

so long as $a^3\not=1$ and $c\not=0$.
Calculating the Jacobian relations we determine that all partial derivatives are $0$ except for
\begin{align*}
-\ti{5}\dP{1}+\dQ{1}=-\frac{(-1+a^3)^3}{4a^4}\text{ and }
-\ti{5}\dP{-1}+\dQ{-1}=\frac{(-1+a^3)^3}{4a^5},
\end{align*}
which are clearly not always $0$.  

\subsection{Determining $Q$}\label{detq}  
For the proofs of Lemma $\ref{ideallemma}$ and subsequently Theorem $\ref{ranktwo}$ to be complete, it only remains to establish that there exists $Q\in R$ so $Q-t_{(5)}t_{(-5)}\in 
\ker(\Pi)$.

Before doing so, we state and prove the following technical fact, which may be found in \cite{N}.

\begin{fact}\label{fact}
Define a bilinear form on the vector space of $n\times n$ matrices over $\C$
by $$\mathbb{B}(A,B)=n\tr{AB}-\tr{A}\tr{B}.$$
Then given vectors $A_1,...,A_{n^2},B_1,...,B_{n^2},$ the $n^2\!\times\! n^2$ matrix $\mathbb{\Lambda}=\bigg(\mathbb{B}(A_i,B_j)\bigg)$
is singular. 
\end{fact}   
\begin{proof}
Consider the co-vector $$v(\ \ )=\left[\begin{array}{c}\mathbb{B}(A_1,\ \ )\\\mathbb{B}(A_2,\ \ )\\
\vdots\\\mathbb{B}(A_{n^2},\ \ )\end{array}\right].$$  If $B_1,...,B_{n^2}$ are linearly dependent then so are
$v(B_1),v(B_2),...,v(B_{n^2})$, which implies the
columns of $\mathbb{\Lambda}$ are linearly dependent.  Otherwise there exists coefficients, not all zero, so
$$c_1 B_1+c_2 B_2+\cdots+c_{n^2}B_{n^2}=\id,$$ which implies $$c_1 v(B_1)+c_2 v(B_2)+\cdots+c_{n^2}v(B_{n^2})=0$$ since the identity $\id$ is
in the kernel of $\mathbb{B}(A,\  )$.  So again the columns of $\mathbb{\Lambda}$ are linearly dependent.  Either way, $\mathbb{\Lambda}$ is
singular.  
\end{proof}

\begin{lem}\label{qlemma} There exists a polynomial $Q\in R$ so $Q-t_{(5)}t_{(-5)}\in \ker(\Pi)$, and in particular
\begin{align}\label{q}
Q=&9-6t_{(1)}t_{(-1)} -6t_{(2)}t_{(-2)} -6t_{(3)}t_{(-3)} -6t_{(4)}t_{(-4)}+t_{(1)}^3 +t_{(2)}^3 +t_{(3)}^3\nonumber\\ 
  &+t_{(4)}^3 +t_{(-1)}^3 + t_{(-2)}^3 +t_{(-3)}^3 +t_{(-4)}^3 -3t_{(-4)}t_{(-3)}t_{(-1)} -3t_{(4)}t_{(3)}t_{(1)}\nonumber\\ 
  &-3t_{(-4)}t_{(2)}t_{(3)} -3t_{(4)}t_{(-2)}t_{(-3)}+3t_{(-4)}t_{(-2)}t_{(1)} +3t_{(4)}t_{(2)}t_{(-1)}\nonumber\\
  &+ 3t_{(1)}t_{(2)}t_{(-3)} +3t_{(-1)}t_{(-2)}t_{(3)}+t_{(-2)}t_{(-1)}t_{(2)}t_{(1)}+t_{(-3)}t_{(-2)}t_{(3)}t_{(2)}\nonumber\\ 
  &+t_{(-4)}t_{(-1)}t_{(4)}t_{(1)} +t_{(-4)}t_{(-2)}t_{(4)}t_{(2)} +t_{(-3)}t_{(-1)}t_{(3)}t_{(1)}+t_{(-3)}t_{(-4)}t_{(3)}t_{(4)}\nonumber\\
  &+t_{(-4)}^2t_{(-3)}t_{(-2)}  +t_{(4)}^2t_{(3)}t_{(2)} +t_{(-1)}^2t_{(-2)}t_{(-4)} +t_{(1)}^2t_{(2)}t_{(4)}+t_{(1)}t_{(-2)}^2t_{(-3)}\nonumber\\ 
  &+t_{(-1)}t_{(2)}^2t_{(3)} +t_{(-4)}t_{(-3)}t_{(1)}^2 +t_{(4)}t_{(3)}t_{(-1)}^2 +t_{(-4)}t_{(2)}t_{(-3)}^2 +t_{(4)}t_{(-2)}t_{(3)}^2 \nonumber\\ 
  &+t_{(-1)}^2t_{(-3)}t_{(2)} +t_{(1)}^2t_{(3)}t_{(-2)} +t_{(-4)}t_{(1)}t_{(2)}^2 +t_{(4)}t_{(-1)}t_{(-2)}^2+t_{(-4)}t_{(3)}t_{(-2)}^2\nonumber\\
  &+t_{(4)}t_{(-3)}t_{(2)}^2 +t_{(1)}t_{(3)}t_{(-4)}^2 +t_{(-1)}t_{(-3)}t_{(4)}^2 +t_{(-1)}t_{(-4)}t_{(3)}^2\nonumber\\
  &+t_{(1)}t_{(4)}t_{(-3)}^2-2t_{(-3)}^2t_{(-2)}t_{(-1)}-2t_{(3)}^2t_{(2)}t_{(1)} -2t_{(-4)}^2t_{(-1)}t_{(2)}\\ 
  &-2t_{(4)}^2t_{(1)}t_{(-2)}+t_{(-1)}^2t_{(-2)}^2t_{(-3)}+t_{(1)}^2t_{(2)}^2t_{(3)}+t_{(-4)}t_{(-1)}^2t_{(2)}^2\nonumber\\
  &+t_{(4)}t_{(1)}^2t_{(-2)}^2-t_{(-4)}t_{(-2)}^2t_{(2)}t_{(1)} -t_{(4)}t_{(2)}^2t_{(-2)}t_{(-1)}\nonumber\\
  &-t_{(-3)}t_{(1)}^2t_{(-1)}t_{(2)}-t_{(3)}t_{(-1)}^2t_{(1)}t_{(-2)}-t_{(-3)}t_{(2)}^2t_{(-2)}t_{(1)} -t_{(3)}t_{(-2)}^2t_{(2)}t_{(-1)}\nonumber\\
  &-t_{(-4)}t_{(-2)}t_{(-1)}t_{(1)}^2 -t_{(4)}t_{(2)}t_{(1)}t_{(-1)}^2-t_{(-1)}t_{(-2)}^3t_{(1)}-t_{(-1)}t_{(2)}^3t_{(1)} \nonumber\\
  &-t_{(-1)}^3t_{(-2)}t_{(2)}-t_{(1)}^3t_{(-2)}t_{(2)}-t_{(-4)}t_{(-3)}t_{(-2)}t_{(-1)}t_{(2)}-t_{(4)}t_{(3)}t_{(2)}t_{(1)}t_{(-2)}\nonumber\\
 &-t_{(-1)}t_{(1)}t_{(2)}t_{(-4)}t_{(3)} -t_{(-1)}t_{(1)}t_{(-2)}t_{(4)}t_{(-3)}+ t_{(-2)}t_{(-1)}^2t_{(1)}^2t_{(2)} +t_{(-1)}t_{(-2)}^2t_{(2)}^2t_{(1)}.\nonumber
\end{align}
\end{lem}

\begin{proof}
The following argument is an adaptation of an existence argument given in \cite{N}, which we use not only to show existence of $Q$, but to derive the explicit 
formulation of $Q$ as well.  Indeed, with respect to Fact $\ref{fact}$, let 
\begin{align*}
&A_1=B_1=\xb_1\ \ \ A_4=B_4=\xb_2^{-1} \ \ \ A_7=B_7=\xb_1\xb_2^{-1}\\
&A_2=B_2=\xb_2 \ \ \ A_5=B_5=\xb_1\xb_2 \ \ A_8=B_8=\xb_2^{-1}\xb_1\\
&A_3=B_3=\xb_1^{-1} \ \ A_6=B_6=\xb_2\xb_1 \ \ A_9=B_9=\xb_2\xb_1^{-1}.
\end{align*}

Then we see that $\mathbb{\Lambda}$ has exactly two entries with 
$\tr{\mathbf{x}_1\mathbf{x}_2\mathbf{x}_1^{-1}\mathbf{x}_2^{-1}}$. After rewriting all matrix entries in terms of our 
generators of $\C[\X]$, we have
$$0 =\det(\mathbb{\Lambda})=P_1\cdot\tr{\mathbf{x}_1\mathbf{x}_2\mathbf{x}_1^{-1}\mathbf{x}_2^{-1}}^2 
+P_2\cdot\tr{\mathbf{x}_1\mathbf{x}_2\mathbf{x}_1^{-1}\mathbf{x}_2^{-1}}+P_3,$$  
where $P_1,P_2,P_3$ are polynomials in terms of 
$$\tilde{R}=\{\tr{\mathbf{x}_1}, \tr{\mathbf{x}_1^{-1}}, 
\tr{\mathbf{x}_2}, \tr{\mathbf{x}_2^{-1}}, \tr{\mathbf{x}_1\mathbf{x}_2}, 
\tr{\mathbf{x}_1^{-1}\mathbf{x}_2^{-1}},\tr{\mathbf{x}_1\mathbf{x}_2^{-1}}, 
\tr{\mathbf{x}_1^{-1}\mathbf{x}_2}\}.$$

If $P_1=0$ then we have a non-trivial relation among the elements of $\tilde{R}$, which we have already seen cannot exist.  
Alternatively, one can evaluate the elements of $\tilde{R}$ with the aid of a computer algebra system to verify that 
$P_1\not=0$.  Then by direct calculation, using {\it Mathematica}, we verify that $P_2=-P \cdot P_1$.  Hence it follows that 
$$-P_3=P_1(t_{(5)}^2-Pt_{(5)})=P_1(t_{(5)}^2-(t_{(5)}+t_{(-5)})t_{(5)})=-P_1t_{(5)}t_{(-5)},$$ and so we have shown the existence 
of $$Q=t_{(5)}t_{(-5)}.$$  Lastly, we simplify $P_3/P_1,$ with the aid of {\it 
Mathematica}, which turns out to be equation $\eqref{q}$.
\end{proof}

\section{Outer Automorphisms}\label{outer}
Given any $\alpha\in \mathrm{Aut}(\F_2)$, we define $a_\alpha \in \mathrm{End}(\C[\X])$ by extending the following mapping 
$$a_\alpha(\tr{\mathbf{w}})=\tr{\alpha(\mathbf{w})}.$$  If $\alpha \in \mathrm{Inn}(\F_2)$, then there exists $\mathtt{u}\in \F_2$ so for all 
$\mathtt{w}\in \F_2$, $$\alpha(\mathtt{w})=\mathtt{uwu}^{-1},$$ which implies 
$$a_\alpha(\tr{\mathbf{w})}=\tr{\ub\mathbf{w}\ub^{-1}}=\tr{\mathbf{w}}.$$

Thus $\ot(\F_2)$ acts on $\C[\X]$.  By results of Nielsen (see \cite{MKS}, \cite{Ni}), $\ot(\F_2)$ is generated by the following mappings
\begin{align}
\tau&=\left\{
\begin{array}{l}
\xt_1\mapsto \xt_2\\
\xt_2\mapsto \xt_1 
\end{array}\right.\\
\iota&=\left\{
\begin{array}{l} 
\xt_1\mapsto \xt_1^{-1}\\
\xt_2\mapsto \xt_2          
\end{array}\right.\\
\eta&=\left\{
\begin{array}{l} 
\xt_1\mapsto \xt_1\xt_2\\
\xt_2\mapsto \xt_2          
\end{array}\right.
\end{align}

Let $\mathfrak{D}$ be the subgroup generated by $\tau$ and $\iota$, and let $\C\mathfrak{D}$ be the corresponding group ring.  
Then $\C[\X]$ is a $\C\mathfrak{D}$-module.

\begin{lem}\label{preservelemma}
The action of $\C \mathfrak{D}$ preserves $R$, and $\mathfrak{D}$ fixes $P$ and $Q$. 
\end{lem}
\begin{proof}
First we note that it suffices to check $$\{\iota, \tau\}$$ on $$\{t_{(\pm i)},\  1\leq i\leq 4\},$$ since the former  
generates $\C\mathfrak{D}$ and the latter generates $R$. Secondly we observe that both $\iota$ and $\tau$ are 
idempotent.

Indeed, $\iota$ maps the generators of $R$ as follows:
\begin{eqnarray*}
t_{(1)}\mapsto t_{(-1)} \mapsto t_{(1)}\\  
t_{(3)}\mapsto t_{(-4)}\mapsto t_{(3)}\\
t_{(2)}\mapsto t_{(2)}\\  
t_{(-2)}\mapsto t_{(-2)}\\
t_{(4)}\mapsto t_{(-3)}\mapsto t_{(4)}.
\end{eqnarray*}
Likewise, $\tau$ maps the generators of $R$ by:
\begin{eqnarray*}
t_{(1)}\mapsto t_{(2)} \mapsto t_{(1)}\\
t_{(-1)}\mapsto t_{(-2)}\mapsto t_{(-1)}\\
t_{(3)}\mapsto t_{(3)}\\
t_{(-3)}\mapsto t_{(-3)}\\
t_{(4)}\mapsto t_{(-4)}\mapsto t_{(4)}. 
\end{eqnarray*}
Hence both map into $R$.

For the second part of the lemma, it suffices to observe $\iota(t_{(\pm 5)})=t_{(\mp 5)}=\tau(t_{(\pm 5)}),$ because in $\C[\X]$, $$P=t_{(-5)}+t_{(5)}\ \text{and}\ 
Q=t_{(5)}t_{(-5)}.$$
\end{proof}

Observing $\iota(t_{(5)})=\tau(t_{(5)})=t_{(-5)}=P-t_{(5)},$ it is apparent that $\mathfrak{D}$ does not act as a permutation group on the
entire coordinate ring of $\X$.  However, when restricted to $R$ there is  

\begin{thm}\label{symmetrytheorem}
$\mathfrak{D}$ restricted to $R$ is group isomorphic to the dihedral group, $D_4$, of order $8$.  Moreover, the algebraically independent 
generators are characterized as those which $\mathfrak{D}$ acts on as a permutation group.
\end{thm}

\begin{proof}
Let $S=\mathrm{Sym}(\pm 1,\pm 2, \pm 3, \pm 4)$ be the symmetric group of all permutations on the eight letters $\pm i$ for $1\leq i\leq 4$.  
Then we have worked out, in the proof of Lemma \ref{preservelemma}, that $\tau$ acts on the subscripts of $t_{(\pm i)}$ as the 
permutation 
$$(1,2)(-1,-2)(4,-4)$$ and likewise, $\iota$ acts as the permutation $$(1,-1)(3,-4)(-3,4).$$

Since $\mathfrak{D}$ is generated by these elements, we certainly have a well defined injection $\mathfrak{D}\to S$.  The Cayley 
table for $\mathfrak{D}$ is:

\begin{center}
{\tiny
\begin{tabular}{|c||c|c|c|c|c|c|c|c|}
\hline
&$id$ & $\iota$ & $\tau$ & $\iota\tau$ & $\tau\iota$ & $\tau\iota\tau$ & $\iota\tau\iota$ & $\tau\iota\tau\iota$\\
\hline
\hline
$id$ & $id$& $ \iota$ & $\tau$ &  $\iota \tau$ & $\tau \iota$ & $\tau \iota \tau$ & $\iota\tau\iota$ & $\tau\iota\tau\iota$ \\
\hline
$\iota$ &$\iota$ & $id$ & $\iota\tau$ & $\tau$ & $\iota\tau\iota$ & $\tau\iota\tau\iota$ &  $\tau\iota$  & $\tau \iota \tau$ \\
\hline
$\tau$ & $\tau$ & $\tau\iota$ & $id$ & $\tau\iota\tau$ & $\iota$ & $\iota\tau$ & $\tau\iota\tau\iota$ & $\iota \tau \iota$\\
\hline
$\iota\tau$ & $\iota\tau$ & $\iota\tau\iota$ & $\iota$ & $\tau\iota\tau\iota$ & $id$ & $\tau$ & $\tau\iota\tau$ &$\tau\iota$ \\
\hline
$\tau \iota$ &$\tau \iota$ & $\tau$ &$\tau \iota \tau$ & $id$ & $\tau\iota\tau\iota$ &$\iota\tau\iota$ & $\iota$ &$\iota\tau$\\
\hline
$\tau \iota \tau$ & $\tau \iota \tau$ &$\tau\iota\tau\iota$ &$\tau \iota$ &$\iota\tau\iota$ &$\tau$ & $id$ &$\iota\tau$ & $\iota$\\
\hline
$\iota \tau \iota$ & $\iota \tau \iota$ &$\iota\tau$ &$\tau\iota\tau\iota$ &$\iota$ & $\tau\iota\tau$ & $\tau\iota$ & $id$ & $\tau$\\
\hline
$\tau\iota\tau\iota$ &$\tau\iota\tau\iota$ & $\tau\iota\tau$ & $\iota\tau\iota$ & $\tau\iota$ &$\iota\tau$ & $\iota$ & $\tau$ &$id$\\
\hline
\end{tabular}}  
\end{center}
where

\begin{center}
\begin{tabular}{ll}
$id\mapsto (1)$ & $\iota\mapsto (1,-1)(3,-4)(-3,4)$\\
$\tau\mapsto (1,2)(-1,-2)(4,-4)$ & $\iota \tau\mapsto (1,2,-1,-2)(3,-4,-3,4)$\\
$\tau \iota\mapsto (1,-2,-1,2)(3,4,-3,-4)$ & $\tau \iota \tau\mapsto (2,-2)(3,4)(-3,-4)$\\
$\iota\tau\iota\mapsto (1,-2)(2,-1)(3,-3)$ & $\tau\iota\tau\iota\mapsto (1,-1)(2,-2)(3,-3)(4,-4)$.
\end{tabular}
\end{center}

It is an elementary exercise in group theory (see \cite{H}) to show any group presentable as $$\{a,b\ \big|\ |a|=n\geq 3,\ 
|b|=2,\ ba=a^{-1}b\}$$ is isomorphic to the dihedral group $D_n$ of order $2n$.  However, letting $a=\tau\iota$ and $b=\iota$ we 
see $|a|=4$, $|b|=2$, $\frak{D}$ is generated by $a$ and $b$, and  $$ba=\iota\tau\iota=(\tau\iota)^{-1}\iota=a^{-1}b.$$    

The last statement in the theorem follows from the fact that $\{t_{(\pm i)}\ |\ 1\leq i\leq 4\}$ are algebraically independent and $\mathfrak{D}$ does not act as a permutation 
group if $t_{(5)}$ were included.
\end{proof}

\begin{rem}
The action of $\frak{D}$ on $\C[\X]$ determines an action on $\X$.  Since $\frak{D}$ acts as a permutation group on $R$ the surjection from Theorem $\ref{ranktwo}$, $\X\to \C^8$, is 
$\frak{D}$-equivariant.  In this way $\X$ exhibits $8$-fold symmetry.
\end{rem}

As already noted, the group ring $\C\mathfrak{D}$ acts on $\C[\X]$.  By brute force computation, one can establish the following succinct 
expressions for the polynomial relations $P$ and $Q$.

\begin{cor}
In $\C\mathfrak{D}$ define $\mathbb{S}_\mathfrak{D}$ to be the group ``symmetrizer'' $$\sum_{\sigma\in \mathfrak{D}}\sigma.$$  Then  
$P=\mathbb{S}_\mathfrak{D}(p)-3$ and $Q=\mathbb{S}_\mathfrak{D}(q)+9$ where $p$ and $q$ are given by:
\begin{align*}
p=&\frac{1}{8}\left(t_{(1)}t_{(-1)}t_{(2)}t_{(-2)}-4t_{(1)}t_{(-2)}t_{(-4)}+2t_{(1)}t_{(-1)}+2t_{(3)}t_{(-3)}\right)\\
q=&\frac{1}{8}\big(2t_{(-2)}t_{(-1)}^2t_{(1)}^2t_{(2)}+4t_{(1)}^2t_{(2)}^2t_{(3)}-4t_{(1)}^3t_{(-2)}t_{(2)}-8t_{(-4)}t_{(-2)}t_{(-1)}t_{(1)}^2-4t_{(4)}t_{(3)}t_{(2)}t_{(1)}t_{(-2)}+\\ 
&8t_{(1)}t_{(3)}t_{(-4)}^2+8t_{(-4)}t_{(1)}t_{(2)}^2 -8t_{(3)}^2t_{(2)}t_{(1)} +4t_{(4)}t_{(-3)}t_{(2)}^2+t_{(-2)}t_{(-1)}t_{(2)}t_{(1)}+t_{(-3)}t_{(-4)}t_{(3)}t_{(4)} + \\
&4t_{(-3)}t_{(-1)}t_{(3)}t_{(1)} +4t_{(1)}^3 +4t_{(3)}^3+12t_{(-4)}t_{(-2)}t_{(1)}-12t_{(-4)}t_{(2)}t_{(3)}-12t_{(1)}t_{(-1)} -12t_{(3)}t_{(-3)}\big).
\end{align*}
\end{cor}
\begin{proof}
We work out $P$ only since the computation for $Q$ is established in the same way but longer.  Indeed,
\begin{align*}
\mathbb{S}_\mathfrak{D}(p)=&\frac{1}{8}\big(\mathbb{S}_\mathfrak{D}(t_{(1)}t_{(-1)}t_{(2)}t_{(-2)})
-4\mathbb{S}_\mathfrak{D}(t_{(1)}t_{(-2)}t_{(-4)})+2\mathbb{S}_\mathfrak{D}(t_{(1)}t_{(-1)})
+2\mathbb{S}_\mathfrak{D}(t_{(3)}t_{(-3)})\big)\\
=&\frac{1}{8}\big(8t_{(1)}t_{(-1)}t_{(2)}t_{(-2)}-4(2t_{(1)}t_{(2)}t_{(-3)}+2t_{(-1)}t_{(-2)}t_{(3)}\\
&+2t_{(1)}t_{(-2)}t_{(-4)}+2t_{(-1)}t_{(2)}t_{(4)})
+2(4t_{(1)}t_{(-1)}+4t_{(2)}t_{(-2)}+4t_{(3)}t_{(-3)}+4t_{(4)}t_{(-4)})\big)\\
=&P+3.
\end{align*}
With the help of {\it Mathematica} or a tedious hand calculation, the formula for $Q$ is equally verified.
\end{proof}

In \cite{AP} an algorithm is deduced that can be adapted to write {\it minimal} generators for $\C[\X]$ when $\F_r$ is free of arbitrary rank, which we do is an upcoming paper.  It is the hope 
of the author that exploiting symmetries as above will simplify the calculations involved in describing the ideals for free groups of rank $3$ or more.  Consequently, this would allow for 
subsequent advances in determining the defining relations of $\X$ in general.


\begin{thebibliography}{100}

\bibitem{AP} Abeasis, A., and Pittaluga, M., {\it On a minimal set of generators for the invariants of 
$3\times 3$ matrices}, Comm. Alg. $\mathbf{17(2)}$ ($1989$), $487$-$499$

\bibitem{ADS} Aslaksen, H., Drensky, V., and Sadikova, L., {\it Defining relations of invariants of two $3\times 3$ matrices}, 
J. Algebra $\mathbf{298}$ $(2006)$, $41$-$57$

\bibitem{D} Dolgachev, I.,``Lectures on Invariant Theory,'' London Mathematical Lecture Notes Series $296$,
Cambridge University Press, $2003$

\bibitem{DF} Drensky, V., and Formanek, E., ``Polynomial Identity Rings,'' Advanced Courses in Mathematics CRM
Barcelona, Birkh$\ddot{\mathrm{a}}$user Verlag Basel, $2004$

\bibitem{Du} Dubnov, J., Sur une g\'{e}n\'{e}ralisation de l'\'{e}quation de Hamilton-Caley et sur les 
invariants simultan\'{e}s de plusieurs affineurs, Proc. Seminar on Vector and Tensor Analysis, Mechanics 
Research Inst., Moscow State Univ. $\mathbf{2/3}$ ($1935$), $351$-$367$.

\bibitem{E} Eisenbud, D.,  ``Commutative Algebra with a View Toward Algebraic Geometry,''  Graduate Texts in Mathematics No. 
150, Spring-Verlag New York, 1995

\bibitem{H} Hungerford, T., ``Algebra,'' Graduate Texts in Mathematics No. $73$, Spring-Verlag New York, $1974$

\bibitem{L} Lawton, S., {\it $\SL$-Character Varieties and $\mathbb{RP}^2$-Structures on a Trinion}, PhD Thesis, University of Maryland, $2006$

\bibitem{MKS} Magnus, W., Karrass, A., and Solitar, D., ``Combinatorial Group Theory: Presentations of Groups in Terms of 
Generators and Relations,'' Pure and Applied Mathematics Vol. XIII, Interscience Publishers, $1966.$

\bibitem{MS} Marincuk, A. and Sibirskii, K., {\it Minimal polynomial bases of affine invariants of square
matrices of order three}, Mat. Issled. $\mathbf{6}$ ($1971$), $100$-$113$

\bibitem{N} Nakamoto, K., {\it The structure of the invariant ring of two matrices of degree $3$}, J. Pure and
Applied Alg. $\mathbf{166}$ ($2002$), $125$-$148$

\bibitem{Ni} Nielsen, J., {\it Die Isomorphismengruppe der freien Gruppen}, Math. Ann., $\mathbf{91}$ ($1924$), $169$-$209$

\bibitem{PX} Previte, J. and Xia, E., Various letters to the author.

\bibitem{P} Procesi, C., {\it Invariant theory of $n\times n$ matrices}, Advances in Mathematics {\bf $19$}  
($1976$), $306$-$381$

\bibitem{R} Razmyslov, Y., {\it Trace identities of full matrix algebras over a field of characteristic zero.} (Russian)
Izv. Akad. Nauk SSSR Ser. Mat.  $\mathbf{38}$  $(1974)$, $723$-$756$.

\bibitem{Ro} Rowen, L. H., ``Polynomial Identities in Ring Theory,'' Academic Press, Inc., New York, $1980$ 

\bibitem{Sc} Schwartz, R., {\it Spherical CR Geometry and Dehn Surgery}, Ann. of Math. Stud. $\mathbf{165}$, $(2007)$, in Press

\bibitem{S} Shafarevich, I., ``Basic Algebraic Geometry $1$,'' $2^{\text{nd}}$ edition,  Springer-Verlag, Berlin, $1994$

\bibitem{Si} Sikora, A., {\it $SL_n$-Character Varieties as Space of Graphs}, Trans. Amer. Math. Soc. $\mathbf{353}$  $(2001)$,  no. $7$, $2773$-$2804$ 

\bibitem{SR} Spencer,A. and Rivlin, R., {\it Further results in the theory of matrix polynomials}, Arch. Rational Mech. Anal. $\mathbf{4}$ 
($1960$), $214$-$230$

\bibitem{T} Teranishi, Y., {\it The ring of invariants of matrices},  Nagoya Math. J. $\mathbf{104}$  $(1986)$, $149$-$161$.

\bibitem{We} Wen, Z.X., {\it Relations polynomiales entre les traces de produits de matrices}, C.R. Acad. Sci. Paris $\mathbf{318}$ $(1994)$ no. $2$, $99$-$104$

\bibitem{W} Wolfram, S., ``Mathematica: A System for Doing Mathematics by Computer,'' Wolfram Press, $2000$

\end{thebibliography}
\end{document}